\documentclass[12pt,a4paper]{article}
\usepackage{url}
\usepackage[colour,hyperref]{ajr}

\usepackage{defns}

\newcommand{\spde}{\textsc{spde}}

\newcommand{\E}[1]{{\cal H}_{#1}}

\title{A step towards holistic discretisation of stochastic partial
differential equations}

\author{A. J. Roberts\thanks{Department of Mathematics \& Computing, 
University of Southern Queensland, Toowoomba, Queensland 4352,
Australia. \protect\url{mailto:aroberts@usq.edu.au}}}

\begin{document}

\maketitle

\begin{abstract}
The long term aim is to use modern dynamical systems theory to derive
discretisations of noisy, dissipative partial differential equations.
As a first step we here consider a small domain and apply stochastic
centre manifold techniques to derive a model.  The approach
automatically parametrises subgrid scale processes induced by spatially
distributed stochastic noise.  It is important to discretise stochastic
partial differential equations carefully, as we do here, because of the
sometimes subtle effects of noise processes.  In particular we see how
stochastic resonance effectively extracts new noise processes for the
model which in this example helps stabilise the zero solution.
\end{abstract}

\tableofcontents

\section{Introduction}

\begin{figure}
    \centering
    \includegraphics[width=\textwidth]{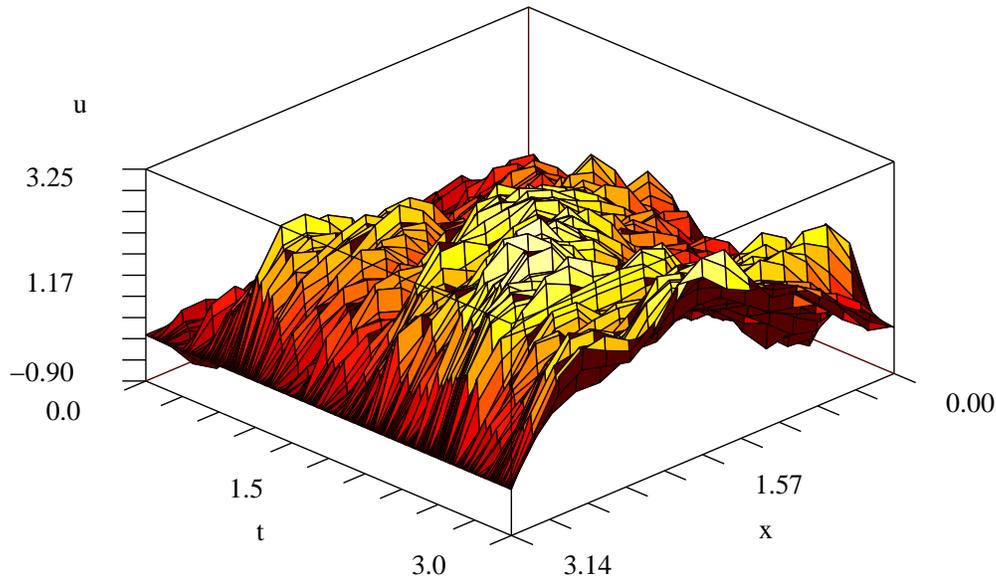}
	\caption{numerical solution over time $0<t<3$ of the
	\spde~(\ref{eq:oburgnm}) on the domain $0<x<\pi$ with stochastic
	forcing~(\ref{eq:onoise}) truncated to the first seven spatial
	modes.  Parameters: $\gamma=0$ so $u\propto\sin x$ is linearly
	neutral although nonlinearly stable; $\sigma=1$ for large forcing;
	numerically $\Delta x=\pi/16$ and $\Delta t=0.01$\,.}
    \label{fig:umesh}
\end{figure}
The ultimate aim is to accurately and efficiently model numerically the
evolution of stochastic partial differential equations (\spde{}s).  An
example solution field~$u(x,t)$, see Figure~\ref{fig:umesh}, shows the
intricate spatio-temporal dynamics typically generated in a \spde.
Numerical methods to integrate stochastic \emph{ordinary} differential
equations are known to be delicate and subtle~\cite[e.g.]{Kloeden92}.
We surely need to take considerable care for \spde{}s as
well~\cite[e.g.]{Grecksch96, Werner97}.

\begin{figure}[tbp]
    \centering
    \includegraphics{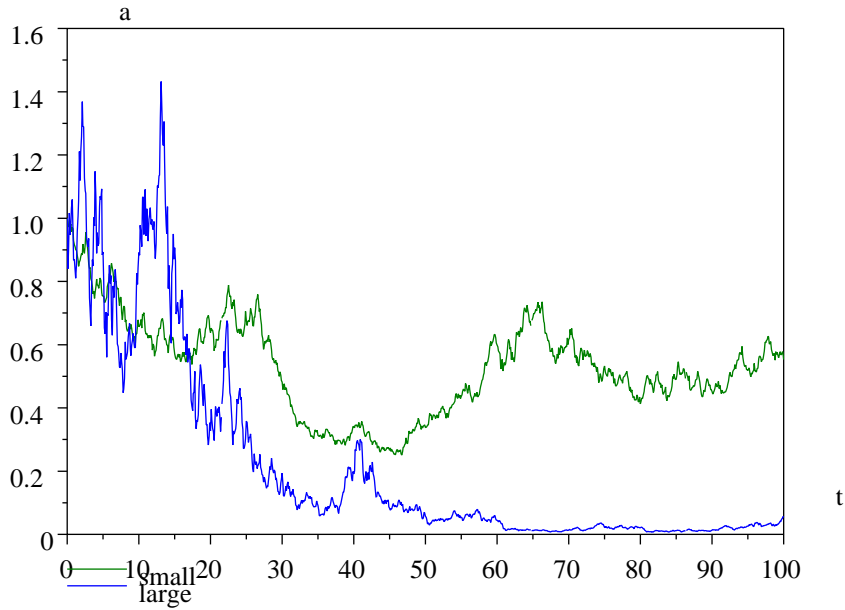}
	\caption{numerical solution of the \sde\ model~(\ref{eq:oomod})
	with small, $\sigma=0.5$, and large, $\sigma=2$, noise.  The
	amplitude~$a$ of the $\sin x$ mode decays for large noise, but not
	for small.  Parameters:  $\gamma=-0.03$ to promote
	linear growth of~$a$, and $\Delta t=0.1$\,.}
    \label{fig:mmodel2}
\end{figure}

An issue is that the stochastic forcing generates high wavenumber,
steep variations, in structures seen in Figure~\ref{fig:umesh}.  Stable
implicit integration in time generally damps far too fast such decaying
modes, yet through stochastic resonance an accurate resolution of the
life-time of these modes may be important on the large scale dynamics.
For example, stochastic resonance causes a high wavenumber noise to
restabilise the trivial solution field~$u=0$ in the simulations
summarised in Figure~\ref{fig:mmodel2}.  Thus we should resolve
reasonably subgrid structures so that numerical discretisation with
large space-time grids achieve efficiency, without sacrificing the
subtle interactions that take place between the subgrid scale
structures.

The methods of centre manifold theory are used here to begin to develop
good methods for the discretisation of \spde{}s.  There is supporting
centre manifold theory by Boxler~\cite{Boxler89, Boxler91, Berglund03}
for the modelling of \sde{}s; the centre manifold approach appears a
better foundation than heuristic arguments for
\sde{}s~\cite[e.g.]{Majda02}.  Further, a centre manifold approach
seems to improve the discretisation of deterministic partial
differential equations~\cite{Roberts98a, Roberts00a, Mackenzie00a,
Roberts01a, Roberts01b, Mackenzie03}.  The first step, taken here, is
to demonstrate the effective modelling of subgrid scale stochastic
structures.

\newpage
\section{Directly seek a one element model}

\begin{figure}[tbp]
    \centering
    \includegraphics{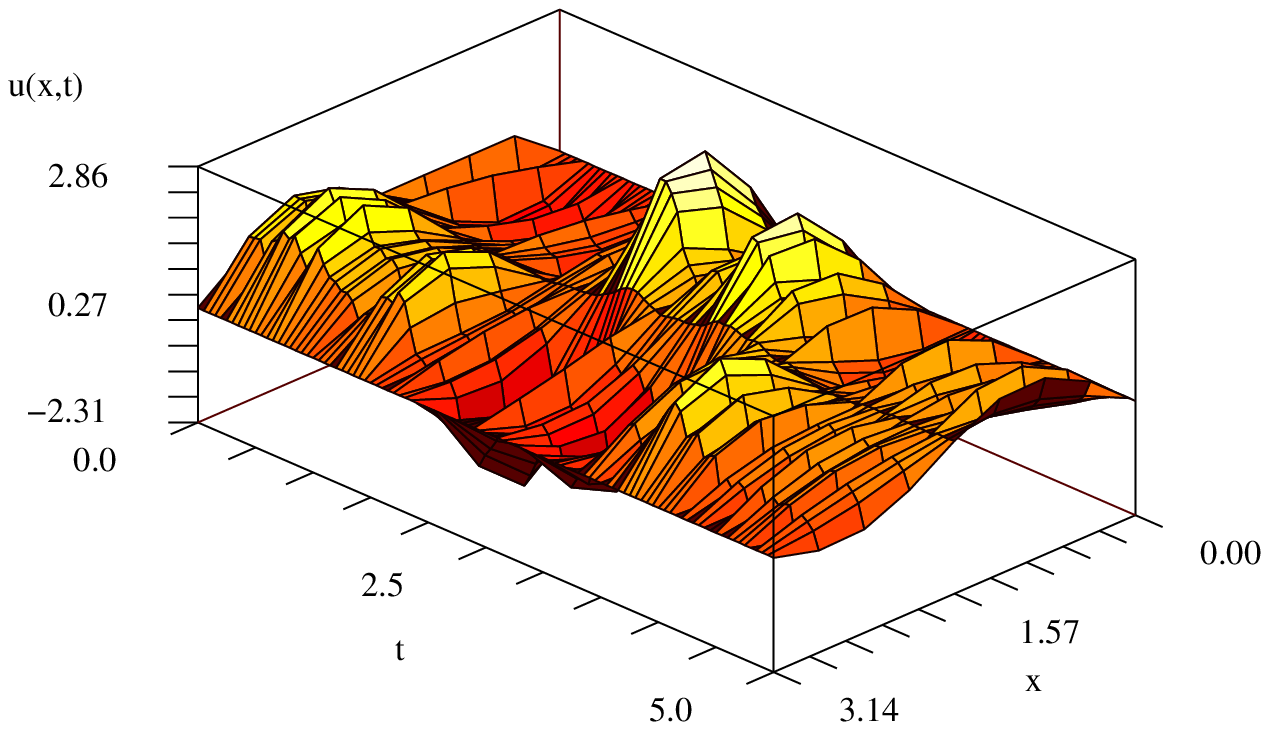}
	\caption{numerical solution of the \spde~(\ref{eq:oburgnm}) with
	strong noise limited to just $\phi=\phi_2(t)\sin 2x$ showing the
	$\sin x$ mode decays.  Parameters: $\sigma=2$, $\gamma=-0.03$ to
	promote linear growth of the $\sin x$ mode, $\Delta t=0.05$ and
	$\Delta x=\pi/8$\,.}
    \label{fig:u1sin2s}
\end{figure}

The simplest case, and that developed here, is the modelling of a
\spde{} on just one finite size element.  Consider the stochastically
forced nonlinear partial differential equation
\begin{equation}
    \D tu=-u\D xu+\DD  x u +(1-\gamma)u +\sigma\phi( x ,t)
    \quad\mbox{such that}\quad  u=0\mbox{ at } x =0,\pi\,,
    \label{eq:oburgnm}
\end{equation}
which involves advection~$uu_x$, diffusion~$u_{xx}$,
reaction~$(1-\gamma)u$, and noise~$\phi$.  In general, the forcing
by~$\phi(x ,t)$, of strength~$\sigma$, is assumed to be white noise
that is delta correlated in both space and time as used in
Figure~\ref{fig:umesh}; however, here we consider only the case
\begin{equation}
    \phi=\phi_2(t)\sin 2x \,,
    \label{eq:onoise}
\end{equation}
where the $\phi_2(t)$ is a white noise that is delta correlated in
time.  Note that the mode $u\propto\sin x$, when $\gamma=0$\,, is
linearly neutral and will form the basis of the model we seek.  Thus
this example of noise forcing the orthogonal $\sin2x$ mode is expected
to be representative of the case of subgrid stochastic forcing and
consequent resolution of higher wavenumber modes.  Many simple
numerical methods, such as Galerkin projection (remembering that the 
domain here represents just one finite element), would completely
obliterate such ``high wavenumber'' modes and hence completely miss
subtle but important subgrid effects.
\begin{figure}[tbp]
    \centering
    \includegraphics{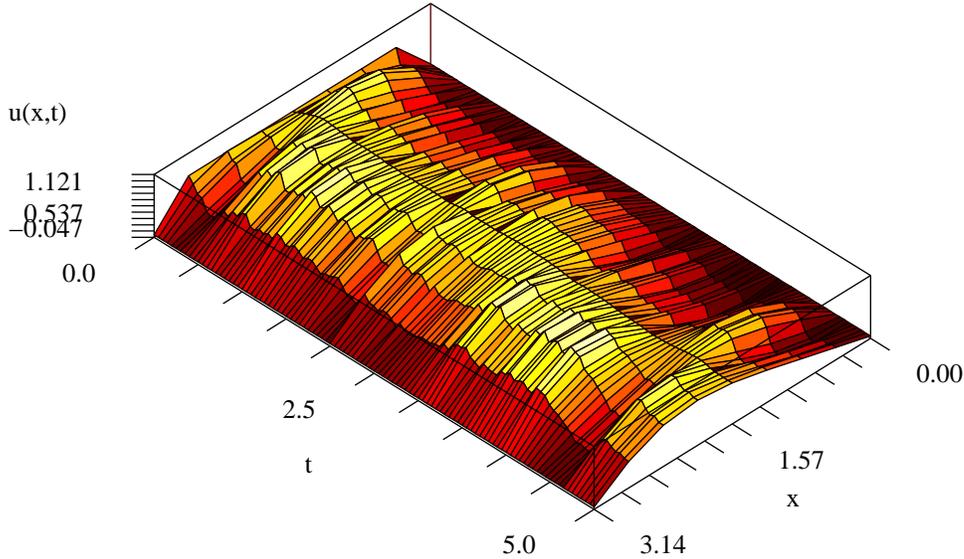}
	\caption{numerical solution of the \spde~(\ref{eq:oburgnm}) with
	relatively weak noise limited to just $\phi=\phi_2(t)\sin 2x$
	showing convergence to a nonlinearly stabilised $\sin x$ mode that
	is perturbed by the noise.  Parameters: $\sigma=0.5$ is small,
	$\gamma=-0.03$ to generate linear growth of the $\sin x$ mode,
	$\Delta t=0.05$ and $\Delta x=\pi/8$\,.}
    \label{fig:u1sin2}
\end{figure}
An example numerical solution, Figure~\ref{fig:u1sin2}, displays that
relatively weak noise only perturbs the deterministic dynamics.
However, when the noise is large enough, then stochastic resonance
restabilises the zero solution and the $\sin x$ mode decays as seen in
Figure~\ref{fig:u1sin2s}.  The success of our approach is seen by it
modelling this induced restabilisation.

For much of the analysis the requirement of white, delta correlated
noise is irrelevant.  Where it is relevant, we interpret the stochastic
differential equations in the Stratonovich sense so that the rules of
traditional calculus apply.

The centre manifold approach identifies that the long term dynamics of
a \spde\ such as~(\ref{eq:oburgnm}) is parametrised by the
amplitude~$a(t)$ of the neutral mode~$\sin x$\,.  Arnold
et~al.~\cite{Arnold95} investigated stochastic Hopf bifurcations this
way, and the approach is equivalent to the slaving principle for
\sde{}s by Schoner and Haken~\cite{Schoner86}.  Computer
algebra~\cite{Roberts96a} determines the solution field
\begin{eqnarray}
    u&=&a\sin x -\rat16a^2\sin2x
    \nonumber\\&&{}
    +\sigma\E2(1-\gamma\E2)\phi_2\sin2x
    -\rat32\sigma a\E3\E2(1-\gamma\E2)\phi_2\sin3x
    \nonumber\\&&{}
    +\rat13\sigma a^2\E4(1+9\E3)\E2\phi_2\sin4x
    +\Ord{a^3+\gamma^2,\sigma^2}\,,
\end{eqnarray}
in which the operator~$\E{m}$ denotes convolution with
$\exp[-(m^2-1)t]$\,.  See in this formula the resolution of the subgrid
structure arising through the interaction of the noise and the
nonlinearity.

The model is the corresponding evolution equation for the amplitude:
\begin{eqnarray}
    \dot a&=& -\gamma a -\rat1{12}a^3
    +\sigma a\rat12\E2(1-\gamma\E2)\phi_2 
    \nonumber\\&&{}
    +\sigma a^3(\rat1{64} +\rat1{12}\E2 -\rat34\E2\E3 +\rat18\E3)\E2\phi_2
    +\Ord{a^4+\gamma^2,\sigma^2}\,.
    \label{eq:onaive}
\end{eqnarray}
This is an unduly messy model as it involves many convolutions over the
rapid time scales we would like to ``step over.''  Straightforward
analyses of forced systems often terminate at this point because of
the tremendously involved form of the repeated convolutions that
occur in higher order terms, especially higher order in the noise
amplitude~$\sigma$.  However, some thought leads us to the drastic
simplifications discussed next.

\newpage
\section{Use a normal form instead}

Here we simplify the model by removing the convolutions from the
evolution equation~(\ref{eq:onaive}).  This step was originally
developed for \sde{}s by Coullet et~al.~\cite{Coullet85} and Sri
Namachchivaya \& Lin~\cite{Srinamachchivaya91}.  In computer algebra
this is done in the equation for the updates to the field and the
evolution:
\begin{displaymath}
    \D t{u'}-\DD x{u'}-u'+{\dot a}'\sin x =\mbox{residual}.
\end{displaymath}
When the residual of the \spde~(\ref{eq:oburgnm}) contains a component
of the form $\E{m}\Phi\sin x $\,, where $\Phi$ denotes some noise
process, which previously we put into ${\dot a}'$ to
form~(\ref{eq:onaive}), we instead recognise that
\begin{equation}
    \frac{d}{dt}\E{m}\Phi =-(m^2-1)\E{m}\Phi+\Phi
    \quad\mbox{thus}\quad
    \E{m}\Phi=\frac{1}{m^2-1}\left[ -\frac{d}{dt}\E{m}\Phi +\Phi
    \right]\,,
    \label{eq:em}
\end{equation}
and so the contribution in the residual is split into: a part
that is integrated into the update~$u'$ for the subgrid field; and a
part without the convolution for the update~${\dot a}'$ for the
evolution.  Note that if the residual component has many convolutions,
then this separation is applied recursively. 

Computer algebra then deduces the normal form model
\begin{equation}
    \dot a=-\gamma a -\rat1{12}a^3
    +\sigma a (\rat16-\rat1{18}\gamma)\phi_2
    -\sigma^2a\rat1{44}(\E2\phi_2-3\E3\E2\phi_2)\phi_2
    +\Ord{a^4+\gamma^2,\sigma^3}\,,
    \label{eq:oomod}
\end{equation}
for the amplitude~$a$ of the $\sin x$ mode, and now to quadratic terms
in the noise.  See that $a=0$ is always a fixed point of this \sde.
Numerical solutions of this model~(\ref{eq:oomod}), see
Figure~\ref{fig:mmodel2}, confirm that for the linearly unstable
(deterministically) parameter~$\gamma=-0.03$ large amounts of noise
restabilise the zero solution.

\section{Stochastic resonance affects deterministic terms} 

The noise~$\phi_2(t)$ so far could have been any distributed forcing at
all, random or deterministic.  The analysis and the results are
generally valid.  We proceed to address the specific modelling
when we restrict the noise~$\phi_2(t)$ to be stochastic white noise in
the Stratonovich sense.  

Previously, the model was a strong model in that~(\ref{eq:oomod}) could
faithfully track given realisations of the original \spde; however, now
we derive the weak model~(\ref{eq:oomodl}) which maintains fidelity to
solutions of the original \spde, but we cannot know which
realisation.

The relevant feature of the large time
model~(\ref{eq:oomod}) is the inescapable and undesirable appearance in
the model of fast time convolutions in the quadratic noise term, namely
$\E2\phi_2 =e^{-3t}\star \phi_2$ and $\E3\E2\phi_2 = e^{-8t}\star
e^{-3t}\star \phi_2$.  These are undesirable because they require
resolution of the fast time response of the system to these fast time
dynamics in order to maintain fidelity with the original
\spde~(\ref{eq:oburgnm}).  However, maintaining fidelity with the full
details of a white noise source is a pyrrhic victory when all we are
interested in is the long term dynamics.  Instead we should only be
interested in those parts of the quadratic noise factors,
$\phi_2\E2\phi_2$ and $\phi_2\E3\E2\phi_2$, that \emph{over long time
scales} are firstly correlated with the other processes that appear and
secondly independent of the other processes: these not only introduce
factors in \emph{new independent} noises into the model but also
introduces a deterministic drift due to stochastic
resonance~\cite[e.g.]{Chao95, Drolet97}.

The argument by Chao \& Roberts~\cite[\S4.1]{Chao95} asserts that we
are interested in the long term statistics of the two quadratic noise
processes $y_1$~and~$y_2$ evolving according to
\begin{equation}
    \dot y_1=z_1\phi_2\,,\quad
    \dot y_2=z_2\phi_2\,,\quad
    \dot z_1=-\beta_1 z_1 +\phi_2\,,\quad
    \dot z_2=-\beta_2 z_2 +z_1\,,
    \label{eq:bin}
\end{equation}
where here the decay rates $\beta_1=3$ and $\beta_2=8$ so that the
convolutions of the noise~$\phi_2$ are represented by the variables
$z_1=\E2\phi_2$ and $z_2=\E3\E2\phi_2$\,.  From the Fokker-Planck
equation for~(\ref{eq:bin}) we have determined that large time solutions
have a probability distribution
\begin{displaymath}
	\mbox{\textsc{pdf}} \propto p(y_1,y_2,t) \exp\left[
	-(\beta_1+\beta_2)z_1^2 +2\beta_2(\beta_1+\beta_2)z_1z_2
	-\beta_2(\beta_1+\beta_2)^2z_2^2 \right]\,,
\end{displaymath}
where the relatively slowly varying~$p$ evolves according to the
approximate equation
\begin{equation}
    \D tp=-\half\D{y_1}p +D:\grad\grad p +\Ord{\grad^3p}
    \label{eq:oofpl}
\end{equation}
where the diffusion matrix
\begin{displaymath}
    D=\left[
    \begin{array}{cc}
        \frac{1}{4\beta_1} & \frac{1}{4\beta_1(\beta_1+\beta_2)} \\
        \frac{1}{4\beta_1(\beta_1+\beta_2)} &
        \frac{1}{4\beta_1\beta_2(\beta_1+\beta_2)}
    \end{array}
    \right]\,.
\end{displaymath}
Interpret~(\ref{eq:oofpl}) as a Fokker-Planck equation and see it
corresponds to the \sde{}s
\begin{equation}
    \dot y_1=\half+\frac{\psi_1(t)}{\sqrt{2\beta_1}}
    \quad\mbox{and}\quad
    \dot y_2=\frac{1}{\beta_1+\beta_2}\left(
    \frac{\psi_1(t)}{\sqrt{2\beta_1}}
    +\frac{\psi_2(t)}{\sqrt{2\beta_2}} \right)\,,
    \label{eq:oosnn}
\end{equation}
where $\psi_i(t)$ are new noises independent of $\phi_2$ \emph{over
long time scales}.  Thus on long time scales, and substituting for
the decay rates~$\beta_i$, we should replace the
quadratic noises by the following:
\begin{equation}
    \phi_2\E2\phi_2=\half+\frac{\psi_1(t)}{\sqrt6}
    \quad\mbox{and}\quad
    \phi_2\E3\E2\phi_2= \frac{\psi_1(t)}{11\sqrt6}
    +\frac{\psi_2(t)}{44} \,.
\end{equation}
Thus the normal form model~(\ref{eq:oomod}) is transformed to
\begin{displaymath}
    \dot a=-\left(\gamma+\rat{\sigma^2}{88}\right) a -\rat1{12}a^3
    +\sigma a (\rat16-\rat1{18}\gamma)\phi_2
    -\sigma^2a(\rat2{121\sqrt6}\psi_1-\rat1{1936}\psi_2) \,.
\end{displaymath}
Combining the new noises into one effective new noise the model is
a little more simply written
\begin{equation}
    \dot a =-\left(\gamma+\rat{\sigma^2}{88}\right) a -\rat1{12}a^3
    +\sigma a (\rat16-\rat1{18}\gamma)\phi_2
    +\sigma^2 a\rat{\sqrt{515}}{1936\sqrt3}\psi \,,
    \label{eq:oomodl}
\end{equation} 
for some white noise~$\psi(t)$ independent of~$\phi_2$ over long times.
Although the nonlinearity induced stochastic resonance generates the
effectively new multiplicative noise, $\propto\sigma^2a\psi$\,, its
most significant effect is the enhancement of the stability of the
equilibrium~$a=0$ through the $\sigma^2a/88$ term.  The equilibrium is
stable for parameters $\gamma>-\sigma^2/88$ which neatly explains the
differences in the stability seen in Figure~\ref{fig:mmodel2} because,
compared to $\gamma=-0.030$, the thresholds for stability are
$-0.003$~and~$-0.045$ for small and large noise respectively.

\begin{figure}[tbp]
    \centering
    \includegraphics{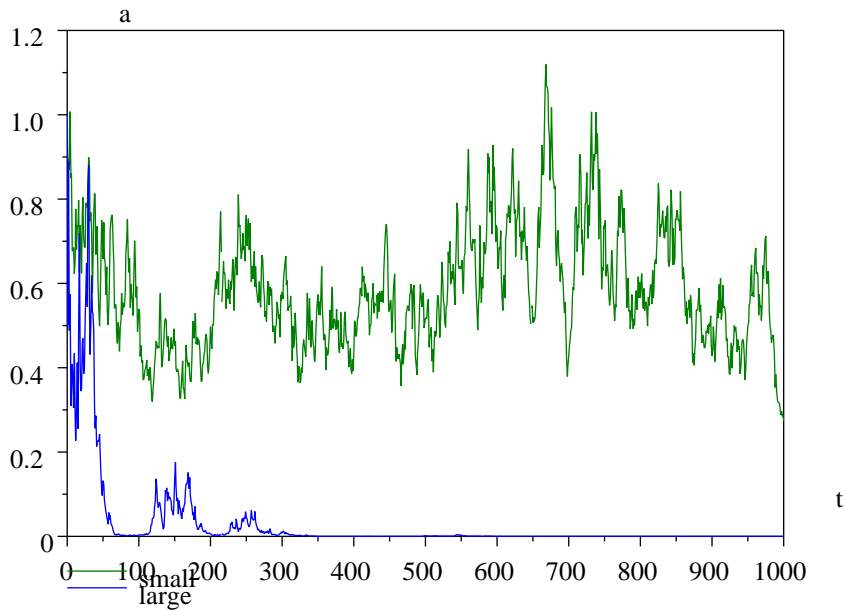}
	\caption{simulations of the long time model~(\ref{eq:oomodl}) for
	small, $\sigma=0.5$, and large, $\sigma=2$, noise over long times.
	Parameters: $\Delta t=1$, $\gamma=-0.03$\,.  }
    \label{fig:mmodelll}
\end{figure}

\section{Conclusion} 

A big virtue of the model~(\ref{eq:oomodl}) is that we may accurately
take large time steps as all the fast dynamics have been eliminated.
Shown in Figure~\ref{fig:mmodelll} are simulations over a long time for
small and large noise again demonstrating the stochastic resonance
induced stabilisation of the equilibrium~$a=0$.  These simulations are
done for an order of magnitude longer times with a time step that is
ten times larger than that we could use previously.

This approach to numerical modelling is viable and effective for
stochastic partial differential equations.  Much more development and
theoretical support is needed.

\bibliographystyle{plain}
\bibliography{ajr,bib,new}

\end{document}